\documentclass[a4paper]{article}
\usepackage{amsmath,amssymb}
\usepackage[all]{xy}
\newtheorem{thm}{Theorem}[section]
\newtheorem{prop}{Proposition}[section]
\newtheorem{cor}{Corollary}[section]
\newtheorem{lem}{Lemma}[section]
\newtheorem{de}{Definition}[section]

\newenvironment{proof}{
                        \noindent{\bf\small Proof.}\small}
                                       {\hfill {$\mathbf \Box$}\medskip}

\newcommand{\got}[1]{\mathfrak{#1}}

\newcommand{\K}{\mathbb{K}}
\newcommand{\N}{\mathbb{N}}
\newcommand{\Z}{\mathbb{Z}}

\newcommand{\Hom}{\mathrm{Hom}}

\newcommand{\g}{\mathfrak{g}}

\makeatletter
\@addtoreset{equation}{section}

\makeatother
\title{On the Generalized Enveloping Algebra of a Color Lie Algebra}
\author{Toukaiddine Petit{{
}}
\\ Freddy Van Oystaeyen\footnote{{\tt Supported by the EC project
 Liegrits MCRTN 505078}.}
\\Departement Wiskunde en Informatica, Universiteit Antwerpen\\
B-2020 (Belgium)\\{toukaiddine.petit@ua.ac.be}\\{fred.vanoystaeyen@ua.ac.be}}
   
\date{}

\begin{document}
\maketitle

\subsubsection*{Abstract}
Let $G$ be an abelien group, $\epsilon$ an anti-bicharacter of $G$ and $L$ a $G$-graded $\varepsilon$ Lie algebra (color Lie algebra) over $\K$ a field of characteristic zero. We prove that all $G$-graded, positive filtered $A$ such that the associated graded algebra  is isomorphic to the $G$-graded $\varepsilon$-symmetric algebra $S(L)$, there is a $G$- graded $\varepsilon$-Lie algebra $L$ and a $G$-graded scalar two cocycle $\omega\in\mathrm{Z}_{gr}^2(L,\K)$, such that $A$ is isomorphic to $ U_\omega(L)$ the generalized enveloping algebra of $L$ associated with $\omega$. We also prove there is an isomorphism of graded spaces between the Hochschild cohomology of the generalized universal enveloping algebra $U(L)$ and the generalized cohomology of color Lie algebra $L$.

\noindent \textbf{Classification AMS 2000}: 16E40, 17B35, 17B75.

\noindent \textbf{Key words}: cohomology, universal enveloping algebra, generalized enveloping algebra, color Lie algebra.
\section*{Introduction}
Let $G$ be an  abelein group, $\epsilon$ an anti-bicharacter of $G$ and  $(L,\left[,\right])$ a $G$-graded $\varepsilon$ Lie algebra (color Lie algebra) over $\K$ a field of characteristic zero. Let $\omega\in\mathrm{Z}_{gr}^2(L,\K)$ be a scalar graded two cocycle of degree zero in the sense of Scheunert-Zhang, \cite{SZ}. The generalized enveloping algebra (or $\omega$-enveloping algebra) of $L$ is the quotient of the $G$-graded tensor algebra $T(L)$  by the $G$-graded two-sided ideal generated by the elements $v_1\otimes v_2-\varepsilon(|v_1|,|v_2|)v_2\otimes v_2-\left[v_1,v_2\right]-\omega(v_1,v_2)$, where $v_1,v_2$ are homogeneous elements of $L$.
The object of the present paper is to study the structure of the generalized enveloping algebra.
In Section 1 we fix notation and
provide background material concerning finite group gradings and color Lie algebras. In Section 2 we introduce the generalized enveloping algebra of color Lie algebra and study its properties. In particular we state the generalized Poincar\'{e}-Birkhoff-Witt theorem for the generalized enveloping algebra. In Section 3 we classify all $G$-graded, positive filtered $A$ such that the associated graded algebra  is isomorphic to the $G$-graded $\varepsilon$-symmetric algebra $S(L)$ which extends the result of Sridharan for Lie algebras, \cite{Sr}. In Section 4 we introduce a graded generalized cohomology (or $\omega$-cohomology) of color Lie algebras which coincides with the graded
Chevalley-Eilenberg cohomology of degree zero of $L$ introduced by Scheunert and
Zhang \cite{SZ} in the case $\omega=0$. We show that there is an isomorphism of graded spaces between the Hochschild cohomology of the generalized universal enveloping algebra and the graded $\omega$-cohomology of color Lie algebra.
\section{Premilinaries}
Throughout this paper  groups are assumed to be abelian and $\K$ is
a field of characteristic zero. We recall some notation for graded
algebras and graded modules \cite{NV}, and some facts on color Lie
algebras from  (\cite{S},\cite{SZ}).
\subsection{Graded Hochschild cohomology}
Let $G$ be a group with identity element $e$. We will write
$G$ as an multiplicative group. An associative algebra $A$ with unit $1_A$, is said to be
$G$-graded, if there is a family $\left\{A_g|  g\in G\right\}$ of
subspaces  of $A$ such that $A=\oplus_{g \in G} A_g$ with  $1_A \in
A_0$ and $A_gA_h \subseteq A_{gh}$, for all $g, h \in G$. Any
element  $a \in A_g$ is  called homogeneous of degree $g$, and we write $|a|=g$.\\
A left graded $A$-module $M$ is a left $A$-module with a
decomposition $M=\oplus_{g \in G}M_g$ such that $A_g.M_h \subseteq
M_{gh}$. Let $M$ and $N$ be graded $A$-modules. Define
\begin{equation}
    \mathrm{Hom}_{A\mbox{-}{\rm gr}}(M,N)=\left\{f\in\mathrm{Hom}_{A}(M,N)|\ \  f(M_g)\subset N_g, \forall\ \ g\in G\right\}.
\end{equation}
We obtain the category of graded left $A$-modules, denoted by $A$-gr,
,\cite{NV}. Denote by ${\rm Ext}_{A\mbox{-} {\rm
gr}}^n(-,-)$ the $n$-th right derived functor of the functor ${\rm
Hom}_{A\mbox{-} {\rm gr}}(-,-)$.\\
Let us recall the notion of graded Hochschild cohomology of a graded algebra $A$.
A graded $A$-bimodule is an $A$-bimodule $M=\oplus_{g \in G}M_g$ such that
 $A_g.M_h.A_k \subseteq M_{ghk}$. Thus we obtain the category of graded $A$-bimodules,
 denoted by $A\mbox{-}A \mbox{-}{\rm gr}$.\\
Let $A^e=A\otimes A^{op}$ be the enveloping algebra of $A$, where
$A^{op}$ is the opposite algebra of $A$. The algebra $A^e$ also is
graded by $G$ by  setting $A^e_g:=\sum_{h \in G} A_h \otimes
A_{h^{-1}g}$. Now the graded $A$-bimodule $M$ becomes a graded
left $A^e$-module by defining the $A^e$-action as
\begin{equation}
    (a \otimes b)m= a.m.b,
\end{equation}
and it is clear that $A^e_g  M_h \subseteq M_{gh}$, i.e., $M$ is a
graded $A^e$-module. Moreover, every graded left $A^e$-module arises
in this way. Precisely, the above correspondence establishes an
equivalence of categories
\begin{align}
    A\mbox{-}A \mbox{-} {\rm gr} \simeq A^e\mbox{-} {\rm gr}.
\end{align}
In the sequel we will identify these categories. Let $M$ be a graded $A$-bimodule, as above, $M$ may be regarded as a graded left
$A^e$-module. The $n$-th graded Hochschild cohomology of $A$ with
value in $M$ is defined by
\begin{align}
\mathrm{HH}_{\rm gr}^n (A, M):= {\rm Ext}_{A^e\mbox{-} {\rm gr}}^n
(A, M), \quad n \geq 0,
\end{align}
where $A$ is the  graded left $A^e$-module induced by the
multiplication of $A$, and the algebra $A^e=\oplus_{g \in G} A^e_g$
is considered as a $G$-graded algebra.
\subsection{Lie color algebras}
The concept of color Lie algebras is related to an abelian group $G$
and an anti-symmetric bicharacter $\varepsilon:G \times G
\rightarrow\K^\times$, i.e.,
 \begin{align}
&\varepsilon\left(g,h\right)\varepsilon\left(h,g\right)=1,\\
    &\varepsilon
\left(g,hk\right)=\varepsilon\left(g,h\right)\varepsilon\left(g,k\right),\\
    &\varepsilon\left(gh,k\right)=\varepsilon\left(g,k\right)\varepsilon\left(h,k\right),
\end{align}
where $g, h, k \in G$ and $\K^\times $ is the multiplicative group
of the units in $\K$.\\
A $G$-graded space $L=\oplus_{g\in G} L_g$ is
said to be a $G$-graded $\varepsilon$-Lie algebra (or simply, color
Lie algebra), if it is endowed with a bilinear bracket
$\left[-,-\right]$ satisfying the following conditions
\begin{equation}
\left[ L_g, L_h \right]\subseteq L_{gh},
\end{equation}
\begin{equation}
    \left[a,b\right]=-\varepsilon\left(|a|,|b|\right)\left[b,a\right],
\end{equation}
\begin{equation}
\varepsilon\left(|c|,
|a|\right)\left[a,\left[b,c\right]\right]+\varepsilon\left(|a|,|b|\right)\left[b,\left[c,a\right]\right]+
\varepsilon\left(|b|, |c|\right)\left[c,\left[a,b\right]\right]=0,
\end{equation}
where $g, h \in G$,  and $a, b, c \in L$ are homogeneous
elements.\\
For example, a
super Lie algebra is exactly  a $\Z_2$-graded $\varepsilon$-Lie
algebra where
\begin{equation}
    \varepsilon(i,j)=(-1)^{ij},\forall\quad i,j\in\Z_2.
\end{equation}
Let $L$ be a color Lie algebra as above and $T(L)$ the tensor
algebra of the $G$-graded vector space $L$. It is well-known that
$T\left(L\right)$ has a natural $\Z\times G$-grading which is fixed
by the condition that the degree of a tensor $a_1\otimes...\otimes
a_n$ with $a_i\in L_{g_i}$, $g_i\in G$, for $1\leq i\leq n$, is
equal to $\left(n,g_1+...+g_n\right)$. The subspace of
$T\left(L\right)$ spanned by  homogeneous tensors of order $\leq n$
will be denoted by $T^n\left(L\right)$. Let $J\left(L\right)$ be the
$G$-graded two-sided ideal of $T\left(L\right)$ which is generated
by
\begin{equation}\label{e1.7}
    a\otimes b-\varepsilon\left(|a|, |b|\right)b\otimes a-\left[a,b\right]
\end{equation}
with homogeneous $a, b \in\g$. The quotient algebra
$U\left(L\right):=T\left(L\right)/J\left(L\right)$ is called the
universal enveloping algebra of the color Lie algebra $L$. The
$\K$-algebra $U\left(L\right)$ is a $G$-graded algebra and  has a
positive filtration by putting $U_n\left(L\right)$ equal to the
canonical image of $T_n\left(L\right)$ in $T\left(L\right)$.\par
 In particular, if
$L$ is $\varepsilon$-commutative (i.e., $[L, L]=0$), then
$U\left(L\right)=S\left(L\right)$ (the $\varepsilon$-symmetric
algebra of the graded  space $L$). \par

The canonical map $\mathrm{i}_L:L\rightarrow U\left(L\right)$ is a
$G$-graded homomorphism and satisfies
\begin{equation}
i_L \left(a\right)i_L\left(b\right)-\varepsilon\left(|a|,
|b|\right)i_L\left(b\right)i_L\left(a\right)=i_L\left(\left[a,b\right]\right).
\end{equation}
The $\Z$-graded algebra $G(L)$ associated with the filtered algebra
$U\left(L\right)$ is defined by letting $G^n\left(L\right)$ be the
vector space $U_n\left(L\right)/U_{n-1}\left(L\right)$ and
$G\left(L\right)$ the space $\oplus_{n\in\N}G^n\left(L\right)$ (note
$U^{-1}\left(L\right):=\left\{0\right\}$). Consequently,
$G\left(L\right)$ is a $\Z\times G$-graded algebra. The well-known
generalized Poincar\'{e}-Birkhoff-Witt theorem, \cite{S}, states
that the canonical homomorphism $\mathrm{i}_L:L\rightarrow
U\left(L\right)$ is an injective $G$-graded homomorphism; moreover,
if $\left\{x_i\right\}_I$ is a homogeneous basis of $L$, where the
index set $I$ well-ordered. Set
$y_{k_j}:=\mathrm{i}\left(x_{k_j}\right)$, then the set of ordered
monomials $y_{k_1}\cdots y_{k_n}$ is a basis of $U\left(L\right)$,
where $k_j\leq k_{j+1}$ and $k_j< k_{j+1}$ if
$\epsilon\left(g_j,g_j\right)\neq 1$ with $x_{k_j}\in L_{g_j}$ for
all $1\leq j\leq n,n\in\N$. In case $L$ is finite-dimensional
$U\left(L\right)$ is a graded two-sided Noetherian algebra (e.g., see for example
\cite{CSV}).

\section{Generalized Enveloping Algebras}
Let $L$ be a $\epsilon$-Lie algebra over $\K$, $U\left(L\right)$ its enveloping algebra and $S\left(L\right)$ its $\epsilon$ symmetric algebra. Let $\omega\in\mathrm{Z}_{gr}^2\left(L,\K\right)$ be a 2-cocycle (of degree zero) for $L$ with values in $\K$ considered as a $G$-graded trivial $L$-module, i.e.
\begin{equation}
\epsilon(|z|,|x|)\omega(x,\left[y,z\right])+\epsilon(|x|,|y|)\omega(y,\left[z,x\right])+\epsilon(|y|,|z|)\omega(z,\left[x,y\right])=0
\end{equation}
for all homogeneous elements $x,y,z\in L$, see \cite{SZ}, \cite{CPV}.
\begin{de}
Let $L$ be a $\epsilon$-Lie algebra and $\omega\in\mathrm{Z}^2_{gr}\left(L,\K\right)$ a scalar graded $2$-cocycle. We call generalized enveloping algebra of $L$ associated with $\omega$, the algebra $U_{\omega}\left(L\right)$, quotient of the tensor algebra over $L$ by the $G$-graded two sided ideal generated by the elements of the form $v_1\otimes v_2-\epsilon(v_1,v_2)v_2\otimes v_1-\left[v_1,v_2\right]-\omega(v_1,v_2)$, where $v_1,v_2$ are homogeneous elements. Then the algebra $U_{\omega}\left(L\right)$ is $G$-graded and $\Z$-filtered.
\end{de}
Let $\omega\in\mathrm{Z}^2\left(L,\K\right)$ be a scalar graded two cocycle of the color Lie algebra $L$. Let $L_{\omega}:=L\ltimes\K\cdot x$ be a central extension of $L$ with $\omega$ such that the new bracket $\left[,\right]'$ is defined by 
\begin{equation} 
   \left[x_1+ax,x_2+bx\right]':=\left[x_1,x_2\right]+\omega\left(x_1,x_2\right)x 
\end{equation}
where $x_1,x_2\in L,a,b,x\in\K$ are homogeneous.
The generalized enveloping algebra $U_{\omega}\left(L\right)$ is isomorphic to the $G$-graded and $\Z$-filtered algebra  $U\left(L_{\omega}\right)/<y-1>$, with $<y-1>$ being the $G$-graded two-sided ideal of $U\left(L_{\omega}\right)$ generated by $y-1$ and $y$ the image of $x$ in $U\left(L_{\omega}\right)$. Denote by 
\begin{equation}
	\pi_{\omega}:U\left(L_{\omega}\right)\stackrel{}{\rightarrow} U_{\omega}\left(L\right)
\end{equation}
the canonical epimorphism.
\begin{de}A graded (left) $(\omega,L)$-module over $\K$ is a graded $\K$-module $M$ endowed with a graded $\K$-linear map $\varphi:L\rightarrow\mathrm{Hom}_{\rm gr}(M,M)$ such that for all homogeneous elements $x,y\in L$
\begin{equation}\label{c2.4}
	\left[[\varphi(x),\varphi(y)]\right]=\varphi(\left[x,y\right])+\omega(x,y)i_M
\end{equation}
where $\left[[\varphi(x),\varphi(y)]\right]=\varphi(x)\varphi(y)-\varepsilon(|x|,|y|)\varphi(y)\varphi(x)$ and $i_M$ is the graded identity map of $M$.
\end{de}
\begin{prop}There is a $1-1$ correspondence between graded (left) $(\omega,L)$-modules and graded (left) $U_\omega(L)$ modules.
\end{prop}
\begin{proof} Let $(M,\varphi)$ be a graded left $(\omega,L)$-module, then the graded $\K$ linear map $\varphi$ may be uniquely extended to a graded $\K$ homomorphism $\widehat{\varphi}:T(L)\rightarrow\mathrm{Hom}_{\rm gr}(M,M)$. It follows from the condition (\ref{c2.4}) that $\widehat{\varphi}$ vanishes on the $G$-graded two sided ideal of $U_\omega(L)$ generated by the elements
$$v_1\otimes v_2-\epsilon(|v|_1,|v|_2)v_2\otimes v_1-\left[v_1,v_2\right]-\omega(v_1,v_2),$$ 
where $v_1,v_2$ are homogeneous elements.
The converse is trivial.
\end{proof}\\
This proves in particular that for any $\omega\in\mathrm{Z}_{gr}^2(L,\K)$ there is a $(\omega,L)$-module, we can see for example that the $\omega$-enveloping algebra $U_\omega(L)$ is a graded $(\omega,L)$-module.
\begin{thm}\label{n1}If $L$ is a $\K$-free $\epsilon$-Lie algebra. Let $\left\{x_i\right\}_{i\in I}$ be a $G$-homogeneous basis of $L$, where $I$ is a well-ordered set. For any central extension of $L$ with $\omega$, the set of ordered monomials $z_{i_1}\cdot\cdot\cdot z_{i_n}$ forms a basis of $U_{\omega}\left(L\right)$, where $i_j\leq i_{j+1}$ and $i_j< i_{j+1}$ if $\epsilon\left(g_j,g_j\right)\neq 1$ with $y_{i_j}\in L_{g_j}$ for all $1\leq j\leq n,n\in\N$.
\end{thm}
\begin{proof}
Since $\left\{x_i\right\}_{i\in I}$ is a $G$-homogeneous basis of the vector space $L$, it follows that $\left\{x_i,x\right\}_{i\in I}$ forms a $G$-homogeneous basis of the vector space $L_{\omega}$. Let
\begin{equation}\label{e2.4}
	i_{\omega}:L_{\omega}\stackrel{i_{L_{\omega}}}{\longrightarrow} U\left(L_{\omega}\right)\stackrel{\pi_{\omega}}{\longrightarrow} U_{\omega}\left(L\right)
\end{equation}
denote the composition. We set $z_i:= i_{\omega}\left(x_i\right)$, $z:= i_{\omega}\left(x\right)$, $y_i:= i_{L_\omega}\left(x_i\right)$, with $i\in I$.
Let $y^{i_0}y_{i_1}\cdots y_{i_n}$  be the generators of the PBW basis of $U\left(L_{\omega}\right)$ with $i_0\in\N$, $i_0\leq i_1$ and $i_0<i_1$ if $\epsilon(|y_{i_0}|,|y_{i_1}|)\neq 1$. In the quotient algebra $U_{\omega}\left(L\right)=U\left(L_{\omega}\right)/<y-1>$, the element $z^{i_0}$ is identified with $1$. Then the canonical projection $\pi_{\omega}$ sends $y^{i_0}y_{i_1}\cdots y_{i_n}$  into $z_{i_1}\cdots z_{i_n}$, and it follows that the elements $z_{i_1}\cdots z_{i_n}$ form a basis of $U_{\omega}\left(L\right)$.
\end{proof}

The restriction of the canonical homomorphism $i_{\omega}$ on $L$, see (\ref{e2.4}), we is again denoted by $i_{\omega}$, i.e., $i_{\omega}:L\rightarrow U_{\omega}\left(L\right)$  satisfies for every $x,y\in L$, homogeneous elements:
\begin{equation}\label{e34}
	 \left[[i_{\omega}(x),i_{\omega}(y)]\right]=i_{\omega}\left(\left[x,y\right]\right)+\omega\left(x,y\right)\cdot i_{U_\omega(L)}
\end{equation}
with $\left[[i_{\omega}(x),i_{\omega}(y)]\right]=i_{\omega}(x)\cdot i_{\omega}(y)-\epsilon\left(|x|,|y|\right)i_{\omega}(y)\cdot i_{\omega}(x)$.

\begin{cor}If $L$ is a $\K$-free $\epsilon$-Lie algebra, then for any central extension of $L$ with $\omega$, $i_{\omega}:L\rightarrow U_{\omega}\left(L\right)$
is an injective homomorphism.
\end{cor}
Thus we may identify every element of $L$ with the canonical image in $U_{\omega}\left(L\right)$. Hence $L$ is embedded in $U_{\omega}\left(L\right)$ and 
\begin{equation}
	\left[[x,y]\right]=\left[x,y\right]+w(x,y)\cdot 1
\end{equation}
for all $x,y\in L$.
The algebra $U_{\omega}\left(L\right)$ has a positive filtration defined by taking for $U_{n,\omega}\left(L\right)$ the canonical image of $U_n(L_{\omega})$ by $\pi_{\omega}$. Denote by $G_\omega(L)$ its associated $\Z$-graded algebra , then $G_{\omega}\left(L\right)$ is a $\Z\times\Gamma$-graded algebra and $\epsilon$-commutative. It follows that the canonical injection 
\begin{equation}
	L\stackrel{i_\omega}{\rightarrow}U_{\omega}(L)\stackrel{}{\rightarrow}G_{\omega}(L),
\end{equation}
may be uniquely extended to a homorphism $\varphi_{\omega}$ of the $\epsilon$-symmetric algebra $S\left(L\right)$ of $L$ into $G_{\omega}\left(L\right)$. If $S^n\left(L\right)$ denotes the set of elements of $S\left(L\right)$ which are homogeneous of degree $\left(n,g_1+...+g_n\right)$, then $\varphi_{\omega}\left(S^n\left(L\right)\right)\subset G_{\omega}^n\left(L\right)$.
\begin{prop}\label{n_2} The canonical homomorphism $\varphi_\omega$ of $S\left(L\right)$ into $G_{\omega}\left(L\right)$ is a $\Z\times\Gamma$-graded algebra isomorphism.
\end{prop}
\begin{proof}
Let $\left\{x_i\right\}_{i\in I}$ be a $G$-homogeneous basis of $L$, with $I$ a well-ordered set. Let $y_{i_1}\cdots y_{i_n}$ be the product $x_{i_1}\cdots x_{i_n}$ calculated in $S\left(L\right)$, $z_{i_1}\cdots z_{i_n}$ the product $x_{i_1}\cdots x_{i_n}$ calculated in $U_{\omega}\left(L\right)$ and $z'_{i_1}\cdots z'_{i_n}$ the canonical image of $z_{i_1}\cdots z_{i_n}$ in $G_{\omega}\left(L\right)$. Since the set of ordered monomials $z_{i_1}\cdot\cdot\cdot z_{i_n}$ form a basis of $U_{\omega}\left(L\right)$, by Theorem \ref{n1}, then the set of ordered monomials $z'_{i_1}\cdot\cdot\cdot z'_{i_n}$ is a basis of $G_{\omega}\left(L\right)$. Since $\varphi_{\omega}\left(y_{i_1}\cdot\cdot\cdot y_{i_n}\right)=z'_{i_1}\cdot\cdot\cdot z'_{i_n}$, it can be seen that $\varphi_{\omega}$ is bijective.
\end{proof}
\begin{prop}If $L$ is of finite dimensional then $U_\omega(L)$ is a graded Noetherian algebra.
\end{prop}
\begin{proof} By Proposition \ref{n_2}, the generalized enveloping algebra $U_\omega(L)$ is a positively graded filtered algebra with its associated graded algebra $gr(U_\omega(L))\simeq S(L)$. The fact that the $\epsilon$-symmetric algebra $S(L)$ is graded Noetherian, see Lemma 2.3 \cite{CSV} and by Theorem 1.1.9 \cite{MB} we deduce that $U_\omega(L)$ is a graded Noetherian algebra.
\end{proof}
\section{Classification of Generalized Enveloping Algebras}
Fix $G$ an abelien group and $\epsilon$ an antisymmetric bicharacter on $G$. Let $V$ be a free $G$-graded vector space over $\K$. Let $S\left(V\right)$ denote the $\epsilon$-symmetric algebra of $V$. Consider the family of all pairs $\left(\mathrm{A},\varphi_{\mathrm{A}}\right)$ where $\mathrm{A}=\cup_{n\in\Z_{+}} F_n\mathrm{A}$ is a $G$-graded, $\Z$-filtered algebra and $\varphi_{\mathrm{A}}:S\left(V\right)\rightarrow G_{F}\left(\mathrm{A}\right)$ is a $G\times\Z$-graded isomorphism. A map $\Psi:\left(A,\varphi_{\mathrm{A}}\right)\rightarrow\left(\mathrm{B},\varphi_{\mathrm{B}}\right)$ is a $G$-graded, $\Z$-filtered algebra homomorphism $\Psi:\mathrm{A}\rightarrow\mathrm{B}$ such that if $G\left(\Psi\right):G\left(\mathrm{A}\right)\rightarrow G\left(\mathrm{B}\right)$ is the $G\times\Z$-graded algebra morphism induced by $\Psi$, the diagram
\begin{equation}
	\xymatrix{
	G(\mathrm{A})
	\ar@{->}[r]^{G(\Psi)} \ar@{<-}[dr]^{\varphi_{\mathrm{A}}}
	& G(B)
		\ar@{<-}[d]^{\varphi_{\mathrm{B}}}\\
	&S(V)}
	\end{equation}
	is commutative. Composition of maps is defined in the obvious way. The resulting category is denoted by $\got{R}_{gr}\left(S(V)\right)$. If $\Psi:\left(A,\varphi_A\right)\rightarrow\left(B,\varphi_B\right)$ is a map then $G\left(\Psi\right):G\left(A\right)\rightarrow G\left(B\right)$ is a graded isomorphism, since $G\left(\Psi\right) =\varphi_{\mathrm{B}}\circ\varphi^{-1}_{\mathrm{A}}$.
\begin{lem}\label{n41}With notation as above $\Psi:\mathrm{A}\rightarrow\mathrm{B}$ is a $\Z$-filtered, $G$-graded isomorphism.
\end{lem}
\begin{proof} Let $\Psi_p:F_p{\mathrm{A}}\rightarrow F_p\mathrm{B}$ denote the $\K$ linear map induced by $\Psi$. We reason by induction on the integer $p$. It is clear that $\Psi_0$ is a graded isomorphism. From the commutativity of the diagram
\begin{equation}
\xymatrix{
0 \ar@{->}[r]&F_{p-1}\mathrm{A}\ar@{->}[r]\ar@{->}[d]^{\Psi_{p-1}}& F_p\mathrm{A}\ar@{->}[r] \ar@{->}[d]^{\Psi_{p}}& G_p(\mathrm{A})\ar@{->}[r]\ar@{->}[d]^{G_p\left(\Psi\right)}& 0\\
0\ar@{->}[r]& F_{p-1}\mathrm{B} \ar@{->}[r]& F_p\mathrm{B}\ar@{->}[r] & G_p(\mathrm{B})\ar@{->}[r]& 0}
\end{equation}
at $\Psi_{p-1}$ and $G_p\left(\Psi\right)$ are graded isomorphisms, it is easily seen that $\Psi_p$ is also a graded isomorphism. Since $p$ is arbitrary, the assertion holds.
\end{proof}
\begin{lem} For each $\left(A,\varphi_A\right)$ pair of $\mathfrak{R}(S(V))$ there is a pair $\left(L,[\omega]\right)$ where $L$ is a $\epsilon$-Lie algebra and $[\omega]\in\mathrm{H}_{gr}^2\left(L,\K\right)$ such that $F_1\mathrm{A}=L_\omega=L\ltimes\K$, with $\omega$ is a representative of $[\omega]$.
\end{lem}
\begin{proof} Let $a,b\in F_1A$ be homogeneous elements, we have $\left[a,b\right]:=ab-\epsilon\left(a,b\right)ba\in F_2 A$. Since $G\left(A\right)$ is $\epsilon$-commutative via $\varphi_A$, then $\left[a,b\right]\in F_1 A$. Thus $F_1 A$ acquires a structure of a $\epsilon$-Lie algebra. It is clear that $\K=F_0 A$ is a central $G$-graded ideal of $F_1A$. The $G$-graded isomorphism $S_1(V)\cong F_1 A/F_0 A$ given by $\varphi_A$, induces a $\epsilon$-Lie structure on $S_1(V)$, denote it by $L$. Then the following sequence 
\begin{equation}
	0\rightarrow\K\stackrel{\mathrm{i}}{\rightarrow}F_1 A\stackrel{\pi}{\rightarrow}L\rightarrow 0
\end{equation}
is central $G$-graded exact and $\pi$ induced by $\varphi_{\mathrm{A}}$. Thus $i$ and $\pi$ are graded homomorphisms  (of degree zero) of $\epsilon$-Lie algebras. Since $S_1 (V)$ is $\K$-free, there exists a graded linear map $\sigma:L\rightarrow F_1 A$ (necessarily of degree zero) such that $\pi\circ\epsilon=\mathrm{id}_{F_1 A}$. We then have
$$\pi([[\sigma(x),\sigma(y)]]-[x,y])=0$$
for all (homogeneous) $x,y\in F_1 A$. Hence, there is a unique map $\omega:L\times L\rightarrow\K$ such that
\begin{equation}\label{e43}
	\mathrm{i}(\omega(x,y))=[[\sigma(x),\sigma(y)]]-\sigma([x,y])
\end{equation}
for all (homogeneous) $x,y\in L$, and it is easy to see that $\omega$ is a homogeneous $2$-cocycle of degree zero, i.e, $\omega\in Z_{gr}^2\left(L,\K\right)$. From \cite{SZ}, it follows that the cohomology class $\left[\omega\right]$ of $\omega$ is independent of the choice of $\omega$.
\end{proof}
\begin{thm}\label{t3.1}Let $G$ be an abelien group and $\epsilon$ a symmetric bicharacter on $G$. Let $V$ be $G$-graded $\K$-free module. Let $S\left(V\right)$ be the $\epsilon$-symmetric algebra on $V$. The isomorphism classes of objects in $\mathfrak{R}_{gr}\left(S(V)\right)$ are in a 1-1 correspondence with pairs $\left(L,\left[\omega\right]\right)$ where $L$ is a $\epsilon$-Lie algebra on $V$ and $\left[\omega\right]$ is an element in $\mathrm{H}_{gr}^2\left(L,\K\right)$. If $\omega$ is a cocycle in the cohomology class $\left[\omega\right]$, then $\left(U_{\omega}\left(L\right),\varphi_{\omega}\right)$ is an object in the isomorphism class determined by $\left(L,\left[\omega\right]\right)$. 
\end{thm}
\begin{proof} Let $L$ be a $\epsilon$-Lie algebra structure on $V$ and $\omega$ is a representative of the cohomology class $\left[\omega\right]\in\mathrm{H}_{gr}^2\left(L,\K\right)$. Using Proposition \ref{n_2}, then $\left(U_{\omega}\left(L\right),\varphi_{\omega}\right)$ is an object in $\mathfrak{R}_{gr}\left(S\right)$. Consider the exact sequence of graded algebras
\begin{equation}
0\rightarrow\K\stackrel{\mathrm{i}}{\rightarrow}F_1\left(U_{\omega}(L)\right)\stackrel{\pi_{\omega}}{\rightarrow}L\rightarrow 0
\end{equation}
where $\pi_{\omega}$ is induced by $\varphi_{\omega}$. The map $i_{\omega}:L\rightarrow F_1\left(U_{\omega}(L)\right)$ is a $\K$-homogeneous linear section and the relation (\ref{e34}) shows that $\left(U_{\omega}\left(L\right),\varphi_{\omega}\right)$ yields $\left(L,\left[\omega\right]\right)$. Let $\left(A,\varphi_A\right)\in\mathfrak{R}_{gr}\left(S(V)\right)$ be another object. Choose $\sigma:L\rightarrow F_1A$ so that (\ref{e43}) is valid for the cocycle $\omega$. Let $\widehat{\sigma}:T\left(L\right)\rightarrow A$ be the natural homogeneous extension of $\sigma$. If $x,y\in L$ are homogeneous, then,
$$ \widehat{\sigma}(x\otimes y-\epsilon(|x|,|y|)y\otimes x-[x,y]-\omega(x,y))=[[\sigma(x),\sigma(y)]]-\sigma([x,y])-\omega(x,y)=0.$$
Then $\widehat{\sigma}$ induces a $G$-graded, $\Z$-filtered homomorphism of algebras $\overline{\sigma}:U_{\omega}(L)\rightarrow\mathrm{A}$. We then have
\begin{equation}
	\xymatrix{
	G(U_{\omega})
	\ar@{->}[r]^{G(\overline{\sigma})} \ar@{<-}[dr]^{\varphi_\omega}
	& G(\mathrm{A})
		\ar@{<-}[d]^{\varphi_{\mathrm{A}}}\\
	&S(V)}
	\end{equation}
For $x\in\g$, $\sigma\left(x\right)$ is in the coset $\varphi_{\mathrm{A}}\left(x\right)$ of $F_1\mathrm{A}$ mod $F_0\mathrm{A}$. Thus,
$G\left(\overline{\sigma}\right)\varphi_\omega\left(x\right)=G\left(\overline{\sigma}\right)i_\omega\left(x\right)=\varphi_{\mathrm{A}}\left(x\right)$. Hence the diagram above is commutative. Thus, $\overline{\sigma}:\left(U_\omega,\varphi_\omega\right)\rightarrow\left(\mathrm{A},\varphi_{\mathrm{A}}\right)$ is a map and then an isomorphism by Lemma \ref{n41}.
\end{proof}\\
From Theorem \ref{t3.1} we retain in particular that $(U_{\omega_1}(L),\varphi_{\omega_1})$ and $(U_{\omega_2}(L),\varphi_{\omega_2})$ are $\Z$-filtered, $G$-graded isomorphic if and only if, $\omega_1$ and $\omega_2$ are (graded) cohomologous. 
 %We give more detail this property by
 %\begin{prop} Let $L$ be a $\K$-free Color Lie algebra and let $\omega_1,\omega_2\in\mathrm{Z}^2_{gr}(L,\K)$ be graded %cohomologous cocycles. Then the maps %\,$\sigma:(U_{\omega_1}(L),\varphi_{\omega_1})\rightarrow(U_{\omega_2}(L),\varphi_{\omega_2})$ are in $1-1$ correspondence with the graded $1$-cochains $f\in\mathrm{C}_{gr}^1(L,\K)$ such that $\omega_2-\omega_1=\delta^1(f)$. The correspondence is defined by the relation
%%\begin{equation}
%	\sigma_f(i_{\omega_1}(x))=f(x)+i_{\omega_2}(x).
%\end{equation}
%Moreover, if $\omega_1,\omega_2,\omega_3\in\mathrm{Z}^2_{gr}(L,\K)$ are such that $\omega_2-\omega_1=\delta^1(f_1)$ and $\omega_3-\omega_2=\delta^1(f_2)$ then $\sigma_{f_1+f_2}=\sigma_{f_1}\circ\sigma_{f_2}$, where $f_1,f_2\in\mathrm{C}_{gr}^1(L,\K)$.
% \end{prop}
\section{Homological Properties of $\mathcal{U}_\omega(L)$ and Color Hopf Algebra}
Let $G$ be a commutative group and $\chi:G\rightarrow\K^*$ a bicharacter. 
 \begin{de}\label{d2} A $(G, \chi)$-Hopf graded algebra $A$ is a 5-tuple $(A, m, \eta, \Delta, \epsilon, S)$ such that  
\begin{enumerate}
	\item $A=\oplus_{g \in G} A_g$ is a graded algebra with multiplication $m : A \otimes A \longrightarrow A$ and
 the unit map $\eta: K \longrightarrow A$. Moreover, $(A, \Delta, \epsilon)$ is a graded coalgebra with
 respect to the same grading.
 \item The counit $\epsilon: A \longrightarrow K$ is an algebra map.
The comultiplication \\$\Delta: A \longrightarrow (A \otimes A)^\chi$ is an algebra map,
where the algebra $(A \otimes A)^\chi$ is equipped
with 
multiplication $*$ defined by 
\begin{align}
(a \otimes b) * (a'\otimes b')= \chi(|b|, |a'|) aa' \otimes bb',
\end{align}
where $a, a' \in A$ and $b, b' \in B$ are  homogeneous.
	\item The antipode  $S: A\longrightarrow A$ is a graded map
such that
 \begin{align}
 \sum a_{1} S(a_{2})=\epsilon(a)=\sum S(a_{1})a_{2}
 \end{align}
 for all homogeneous $a \in A$, where we use Sweedler's notation $$\Delta(a)=\sum a_{1}\otimes a_{2}$$.
\end{enumerate}
\end{de}
\begin{de}An algebra is said to be a color Hopf algebra if it is a $(G,\chi)$-Hopf algebra with the antipode being an isomorphism.
\end{de}
Let $M$ be a graded $A$-bimodule, then we define a left $A$-module by
\begin{align}
 am= \sum \chi (|a_{(2)}|, |m|) a_{(1)}.m. S(a_{(2)}),
\end{align}
for homogeneous $a\in A$ and $m \in M$. It is called the adjoint $A$-graded module and denoted by $^{ad} M$.
\begin{thm}\label{t4}
Let $A=(A, m, \eta, \Delta, \epsilon, S)$ be a color Hopf algebra and let $M$ be a graded $A$-bimodule. Then there exists an isomorphism of graded spaces
\begin{align*}
\mathrm{HH}_{\rm gr}^n(A, M)\simeq{\rm Ext}^n_{A \mbox{-} {\rm gr}}
(\K, ^{ad} M),\quad n\geq 0,
\end{align*}
where $\K$ is viewed as the trivial graded $A$-module via the
counit $\epsilon$, and $^{ad} M$ is the adjoint $A$-module
associated to the graded $A$-bimodule $M$.
\end{thm}
\begin{proof} See \cite{CPV}.
\end{proof}
\begin{prop}Let $L$ be a $\epsilon$-Lie algebra and $\omega\in\mathrm{Z}_{gr}\left(L,\K\right)$ a scalar $2$-cocycle. Then the generalized enveloping algebra $U_\omega(L)$ of $L$ is a color Hopf algebra.
\end{prop}
\begin{proof} It 's shown in \cite{CPV} that the graded tensor algebra $T(L)$ is a color Hopf algebra. Moreover it is easy to prove that the two-sided ideal generated by the elements, $v_1\otimes v_2-\epsilon(v_1,v_2)v_2\otimes v_1-\left[v_1,v_2\right]-\omega(v_1,v_2)$, where $v_1,v_2$ are homogeneous elements, is a graded Hopf ideal. It follows that the generalized enveloping algebra becomes a color algebra Hopf by quotient.
\end{proof}\\
Now we can apply the theorem for the generalized universal enveloping algebra
$U_\omega(L)$ of a $G$-graded $\varepsilon$-Lie algebra $L$.
In fact, if $M$ is a graded  $U_\omega(L)$-bimodule, the adjoint
$U_\omega(L)$-module $^{ad}M$ is given by
\begin{align}\label{e111}
ad(x)m= x.m-\varepsilon(|x|, |m|)m.x=\left[[x,m]\right]
\end{align}
for all homogeneous $x \in L$ and $m \in M$. Thus we have\par
\begin{cor}\label{c}Let $L$ be a $G$-graded $\varepsilon$-Lie algebra and $U_\omega(L)$ its universal generalized enveloping algebra.
Let $M$ be a graded $U_\omega(L)$- bimodule. Then there exists a
graded isomorphism
\begin{align*}
\mathrm{HH}_{\rm gr}^n(U_\omega(L), M)= {\rm Ext}^n_{U_\omega(L)\mbox{-}{\rm gr}}
(\K, {^{ad}M}), \quad n\geq 0,
\end{align*}
where $^{ad} M$ is the adjoint $U_\omega(L)$-module associated with the
graded $U_\omega(L)$-bimodule $M$ defined by (\ref{e111}).
\end{cor}
It follows from \cite{CPV} that the sequence
\begin{equation}
    C:...\rightarrow C_n\stackrel{d_n}{\rightarrow}C_{n-1}\rightarrow...C_1\stackrel{\epsilon}{\rightarrow}C_0
\end{equation}
is a $G$-graded $U(L_\omega)$-free resolution of the $G$-graded trivial $U(L_\omega)$-left module $\K$ via $\epsilon$ where $C_n=U(L_\omega)\otimes_\K\wedge_\varepsilon^n L_\omega$ and the operator $d_n$ is given by
\begin{align*}
&d_n(u\otimes\left\langle x_1, \cdots, x_{n}\right\rangle)\\
&=\sum_{i=1}^{n} (-1)^{i+1} \varepsilon_i\;ux_i\otimes\left\langle x_1, \cdots, \hat{x_i}, \cdots, x_{n}\right\rangle\\
&\quad +\sum_{1\leq i< j \leq n} (-1)^{i+j} \varepsilon_i \varepsilon_j\varepsilon(x_j,x_i)\;u\otimes\left\langle [x_i, x_j], x_1, \cdots, \hat{x_i}, \cdots, \hat{x_j}, \cdots, x_{n}\right\rangle,
\end{align*}
for all homogeneous elements $u \in U(L_\omega)$ and $x_i \in L$, with $\varepsilon_i=\prod_{h=1}^{i-1}\varepsilon(|x_h|, |x_i|)$ $i\geq 2$, $\varepsilon_1=1$ and the sign\quad$\widehat{}$\quad indicates that the element below it must be omitted. The differential operator $d$ maps $U(L_\omega)<y-1>\otimes_\K\wedge_\varepsilon^n L$ into itself, then it passes to the quotient, i.e.,
$\overline{d}:\overline{C}_n\rightarrow \overline{C}_{n-1}$ and satisfies that $\overline{d}\circ \overline{d}=0$ where $\overline{C}_n=U_\omega(L)\otimes_\K\wedge_\varepsilon^n L$.
\begin{prop}The sequence
\begin{equation}
    \overline{C}:...\rightarrow \overline{C}_n\stackrel{\overline{d}_n}{\rightarrow}\overline{C}_{n-1}\rightarrow...\overline{C}_1\stackrel{\epsilon}{\rightarrow}C_0
\end{equation}
is a $G$-graded $U_\omega(L)$-free resolution of the $G$-graded trivial $U_\omega(L)$-left module $\K$ via $\epsilon$.
\end{prop}
\begin{proof} Let $\left\{x_i\right\}_I$ be a homogeneous basis of $L$, where $I$ is a well-ordered set.
By Theorem \ref{n1} the elements
\begin{equation}\label{e113}
    x_{k_1}\cdots x_{k_m}\otimes\left\langle x_{l_1}\cdots x_{l_n}\right\rangle
\end{equation}
with
\begin{equation}
    k_1\leq\cdots\leq k_m\quad\mathrm{and}\quad k_i<k_{i+1}\quad\mathrm{if}\quad\varepsilon(|x_{k_i}|,|x_{k_i}|)=-1
\end{equation}
and
\begin{equation}
    l_1\leq\cdots\leq l_n\quad\mathrm{and}\quad l_i<l_{i+1}\quad\mathrm{if}\quad\varepsilon(|x_{l_i}|,|x_{l_i}|)=1
\end{equation}
form a homogeneous basis of $\overline{C}_n$. The canonical filtration of $U_\omega(L)$, induces a filtration on the complex $\overline{C}$. The associated $Z$-graded complex $G(\overline{C})$ is $G$-graded and isomorphic to the $Z\times G$-graded complex $S(L)\otimes\wedge_\varepsilon L$. It follows from Lemma 3, \cite{CPV}, that the  complex $G(\overline{C})$ is acyclic and consequently so is $\overline{C}$.
\end{proof}

Let $M$ be a $G$-graded left $(\omega,L)$-module, we define the
$n^{th}$ graded cohomology group of $L$ with coefficients in $M$
by
\begin{equation}
     \mathrm{H}_{gr,\omega}^n(L,M):=\mathrm{Ext}_{U_\omega(L)-gr}^n(\K,M).
\end{equation}
The modules on the right hand side can be computed using the left graded $U_\omega(L)$-projective resolution of $\K$. 
If $M$ is a graded left $(\omega,L)$-module, the graded cohomology groups are the graded homology groups of the complex:
$$\Hom_{U_\omega(L)-gr}(\overline{C}_n,M)=\Hom_{U_\omega(L)-gr}(U_\omega(L)\otimes\wedge^n_{\varepsilon}L,M)=\Hom_{\K-gr}(\wedge^n_{\varepsilon}L,M).$$
The coboundary operator in this cocomplex is
\begin{align}
&\overline{\delta}_n(f)\left(x_1, \cdots, x_{n+1}\right)\\
&=\sum_{i=1}^{n+1} (-1)^{i+1}\varepsilon_i\;x_i f\left( x_1, \cdots, \hat{x_i}, \cdots, x_{n+1}\right)\\
&\quad +\sum_{1\leq i< j \leq n+1} (-1)^{i+j}\varepsilon_i \varepsilon_j\varepsilon(x_j,x_i)\;f\left( [x_i, x_j], x_1, \cdots, \hat{x_i}, \cdots, \hat{x_j}, \cdots, x_{n+1}\right).
\end{align}
\begin{thm}\label{t6}
 Let $L$ be a $G$-graded $\epsilon$-Lie algebra, $\omega\in\mathrm{H}_{gr}^2(L,\K)$ and let $U_\omega(L)$ be its generalized universal enveloping algebra. Let $M$  be a graded $U_\omega(L)$-bimodule. Let $^{ad}M$ be the adjoint graded left $(L,\omega)$-module defined by
 $$ad(x)m=\left[[x,m]\right]:=xm-\epsilon(|x|,|m|)mx$$
 for all homogeneous elements $x\in L$ and $m\in M$.
 There exists an isomorphism
 \begin{equation}
 \mathrm{H}_{gr,\omega}^n(L,^{ad}M) \simeq\mathrm{HH}^n_{gr} (U_\omega(L),M), \quad
 n\geq 0.
 \end{equation}
 \end{thm}
 \begin{proof} It is a direct consequence from above and Corollary \ref{c}.
 \end{proof}

\end{document}